\documentclass[12pt]{amsart}
\usepackage{amsfonts}

\theoremstyle{definition}

\theoremstyle{remark}

\newcommand{\const}{\mathop{\rm const}\limits}

\newcommand{\vraisup}{\mathop{\rm vraisup}\limits}

\newcommand{\card}{\mathop{\rm card}\limits}

\newcommand{\degree}{\mathop{\rm degree}\limits}

\begin{document}

\begin{center}

{\bf INTEGRAL OPERATORS IN BILATERAL GRAND LEBESGUE SPACES} \par

\vspace{3mm}
{\bf E. Ostrovsky}\\

e-mail: galo@list.ru \\

\vspace{3mm}

{\bf L. Sirota}\\

e-mail: sirota@zahav.net.il \\

\vspace{3mm}

{\bf E.Rogover} \\

e-mail: rogovee@gmail.com

\vspace{4mm}

 Abstract. \\
{\it In this paper we estimate  the norm of operator acting from one Bilateral Grand
Lebesgue Space (BGLS)  into other Bilateral Grand Lebesgue Space. \par
 We also give some examples to show the sharpness of offered inequalities.} \\

\end{center}

\vspace{3mm}


\vspace{3mm}


\vspace{3mm}

\section{Introduction}

\vspace{3mm}

 Let $  (S_1, \Sigma_1, \mu_1) $  and $  (S_2, \Sigma_2, \mu_2) $ be two measurable
 spaces with sigma-finite non-trivial measures $ \mu_1, \ \mu_2. $ We  denote as usually
 for any measurable function $ f: \Sigma_1 \to R $
 $$
 |f|_{p,\mu_1}= \left[ \int_{S_1} |f(x)|^p \ d\mu_1(x)  \right]^{1/p}, p \in [1,\infty),
 $$
$ f \in L_p(\mu_1) \Leftrightarrow \ |f|_{p,\mu_1} < \infty; $
 $$
 |f|_{p,\mu_2} = \left[ \int_{S_2} |f(x)|^p \ d\mu_2(x)  \right]^{1/p}, p \in [1,\infty),
 $$
or simple

 $$
 |f|_{p,\mu} = \left[ \int_{S} |f(x)|^p \ d\mu(x)  \right]^{1/p}, p \in [1,\infty),
 $$
 $ f \in L_p \Leftrightarrow \ |f|_{p} < \infty $
in the case when $  (S_1, \Sigma_1, \mu_1) =  (S_2, \Sigma_2, \mu_2)
= (S, \Sigma, \mu); $ and  denote $ L(a,b) = \cap_{p \in (a,b)} L_p.  $ \par
 As usually,
 $$
 |f|_{\infty,\mu} = \vraisup_{x \in X} |f(x)|.
 $$
We will denote for simplicity in the case when $ X \subset R^d, d=1,2,\ldots  $
and when $ \mu $ is usually Lebesgue measure

$$
|f|_p = |f|_{p,\mu} = \left[ \int_X |f(x)|^p dx \right]^{1/p}.
$$

 Let $ T $ be linear integral operator of a view:

$$
T[f](s)= Tf(x) = \int_{S_2} K(s, s_2) \ f(s_2) \ d \mu(s_2), \ s \in S_1. \eqno(1)
$$
Here the kernel $ K = K(s_1,s_2), \ s_1 \in S_1, s_2 \in S_2 $ is bimeasurable, i.e.
measurable relative the sigma-algebra $ \Sigma_1 \times \Sigma_2 $ with values in
the real axis $ R $ function. We denote as $ |K|_{p,q}  $ the norm of the operator $ T $
from the space $ L_{p,\mu_2} $ into the space $ L_{q,\mu_1}: $

$$
|T|_{p,q}= \sup_{f \in L_{p,\mu_2}, \ f \ne 0} \frac{ |T[f]|_{q,\mu_1}}{|f|_{p,\mu_2}}.
\eqno(2)
$$

 It is evident that in general case the function $ K = |K|_{p,q} $ may be infinite for
 some values $ (p,q), $ therefore we denote

 $$
  G(T) = \{(p,q): \ p,q \in [1,\infty), |T|_{p,q} < \infty \}
 $$
and suppose $ G(T) \ne \emptyset. $ \par
We define formally for the values $ (p,q) $ which does not belonging the set
$ G(T) \ |T|_{p,q} = + \infty. $  \par
 Denote also for $ j=1,2,3 $
 $$
  G^{(j)}(K) = \{(p,q): \ p,q \in [1,\infty), |K|^{(j)}_{p,q} < \infty \}
 $$
 It is evident that $ G(T) \supset G^{(j)}(K). $  \par

  The complete description of the set $ G(T) $ it follows from the classical interpolation
  Riesz-Thorin theorem, \cite{Adams1}, \cite{Bennet1}.\par
 There are many estimations of the value $ T_{p,q}; $ see, for example,
 \cite{Dunford1}, chapter 5,\cite{Kantorovicz1}, chapter 7-9,  \cite{Okikiolu1},
 chapter 5. Let us denote for simplicity {\it the arbitrary upper estimation} of the value $ |T|_{p,q}  $ as $ K^{(j)}_{p,q}: \ |T|_{p,q} \le K^{(j)}_{p,q}= K_{p,q}. $ \par
  We recall here some expressions for the functionals $ K^{(j)}_{p,q}; $ obviously, these
 formulas are reasonable in the case of finiteness $ K^{(j)}_{p,q}. $ \par

  For instance, $ K^{(1)}_{p,p} = $

 $$
 \left[ \vraisup_{s_2 \in S_2}\int_{S_1}|K(s_1,s_2)|
 \ d \mu(s_1) \right]^{1-1/p} \times
\left[\vraisup_{s_1 \in S_1} \int_{S_2}|K(s_1,s_2)| \ d \mu(s_2) \right]^{1/p};
\eqno(3)
$$

\vspace{3mm}

$$
K^{(2)}_{q,r} = M_1^{1/p_1-p_2/(pp_1)} \cdot  M_2^{1/p}, \eqno(4)
$$
where $ 1 \le p,q, p_1,q_1,p_2,q_2 < \infty, \ pp_1/(p-p_2) = r \in [1,\infty), \ q \le q_2, $

$$
\frac{1}{p_1} + \frac{1}{q_1} =
\frac{1}{p_2} + \frac{1}{q_2}=1 = \frac{1}{p} + \frac{1}{q};
$$

$$
M_1 = \vraisup_{s_1 \in S_1} \int_{S_2} |K(s_1,s_2)|^{p_2} \ d \mu_2(s_2),
$$

$$
M_2 = \vraisup_{s_2 \in S_2} \int_{S_1} |K(s_1,s_2)|^{p_1} \ d \mu_1(s_1);
$$

\vspace{3mm}

$$
K^{(3)}_{p,r} = \left\{ \int_{S_1} d \mu_1(s_1) \cdot
\left[ \int_{S_2} |K(s_1,s_2)|^q \ d \mu_2(s_2) \right]^{r/q} \right\}^{1/r},
\eqno(5)
$$
$ 1/p + 1/q = 1, \ p,q,r \in (1,\infty) $ etc. \par
 Evidently, we can define

 $$
 K^{(4)}_{p,q} = \min_{j=1,2,3} K^{(j)}_{p,q}.
 $$

 Another example (weight Hardy operator).
 Let $ X = (a,b), \ 0 \le a < b \le \infty $ with ordinary Lebesgue
 measure $ d \mu = dx. $ Let us consider the following integral operator, which used
 in the theory of Sobolev spaces \cite{Maz'ya1}:

$$
T_a f(x):= v(x) \int_a^x u(t)f(t) \ dt.
$$
Let $ 1 < q < p \le \infty; \ p^/=p/(p-1), \ q^/ = q/(q-1), \
 \ 1/s:= 1/q \ - \ 1/p, \ B = B(p,q) :=$

$$
\left\{ \int_a^b \left(  \left( \int_x^b |v(t)|^q \ dt \right)^{1/q}
\left(\int_a^x |u(t)|^{p^/} \right)^{1/q^/} \right)^s |u(x)|^{p^/} \ dx
 \right\}^{1/s}.
$$
It is known (see \cite{Edmunds1}) that

$$
q^{1/q}(p^/ q/s)^{1/q^/} \ B \le |T_a|_{p,q} \le q^{1/q} \left(p^/\right)^{1/q^/} \ B.
$$
 In the case when $a=0, b = \infty,  u(t) = 1, v(x) = 1/x $ we obtain the classical
 Hardy operator

 $$
 T_0[ f](x):= x^{-1} \int_0^x f(t) dt
 $$
with the known exact value of $ L_p \to L_p $  norm:
$$
|T_0|_{p,p} = p/(p-1), \ p \in (1,\infty],
$$
 see \cite{Hardy1}, p. 229-238. There are many analogous examples at the same
 place, e.g.:
 $$
 T_+[f](x):= \int_0^{\infty} \frac{f(y) \ dy}{x+y}, \ |T_+|_{p,p} =
 \frac{\pi}{\sin(\pi/p)}, \ p \in (1,\infty) ;
 $$
 $$
 T_m[f](x):= \int_0^{\infty} \frac{f(y) \ dy}{\max(x,y)}, \ |T_m|_{p,p} =
 \frac{p^2}{p-1}, \ p \in (1,\infty),
 $$
 $$
 T_d[f](x) := \int_x^{\infty} \frac{f(y) \ dy}{y}, \ |T_d|_{p,p}=p, \ p \in (1,\infty),
 $$
$$
T_l[f](x) = \int_0^{\infty} \frac{\log(x/y) f(y) \ dy}{ x-y},
|T_l|(p,p) = \left[ \frac{ \pi }{\sin(\pi/p) }  \right]^{2},
$$
  etc. \par

\vspace{3mm}

 Another examples. Let us consider the operators of a "multiplicative"  view:
 $$
 T_M[f](x) = T_Mf(x)= \int_0^{\infty} K(x \cdot y) f(y) dy, \ K(x) \ge 0.
 $$

Denote by

$$
\zeta(s) = \int_0^{\infty} K(x) \ x^{s-1} \ dx
$$
the Mellin's transform of the kernel $ K(\cdot), $ and assume that there exist a
values $ a,b; \ 1 \le a < b \le \infty $ for which
$$
 \forall p \in (a,b) \ \Rightarrow \zeta(1/p) < \infty.
$$
 It is known \cite{Hardy1}, p. 256,  \cite{Mitrinovich1}, p. 146 that

$$
|T_M[f]|^p_p \le \zeta(1/p) \int_0^{\infty} (x |f(x)|)^p \ dx/x^2
$$
and analogously
$$
\int_0^{\infty} x^{p-2} |T_M [f](x)|^p \ dx \le \zeta^p((p-1)/p) \ |f|^p_p.
$$
 In the case when $ K(x) = \exp(-x) $ we obtain the classical Laplace transform:

 $$
 \Lambda[f](x) = \Lambda f(x) = \int_0^{\infty} \exp(-xy) f(y) dy.
 $$
In this case $ \zeta(s) = \Gamma(s) $ and moreover

$$
|\Lambda[f]|_{p^/} \le (2 \pi/p^/)^{1/p^/} \ |f|_p, \ p \in [1,2].
$$
 Analogous estimations are true for Fourier transform. \par
In the books \cite{Hardy1}, chapter 6 and  \cite{Mitrinovich1}, chapter 5 there are
many other examples of integral operators with calculated or estimated norms in the
classical Lebesgue spaces $ L_p.$ \par
  There exists a possibility when for any value of variable $ p, \ p \in (A,B) $ there
 exists (in general case) a {\it unique } value $ q, \ q \in (1,\infty) $ for which
 $  |T|_{p,q} < \infty; $ see example further. We will denote in this case
 $ T = T(\cdot) \in UC $ and will denote the correspondent function by
 $ w: w = w(p): \ W(p):= |T|_{p,w(p)} < \infty  $ and

 $$
 \forall t \ne w(p) \ \Rightarrow |T|_{p,t} = \infty.
 $$

\vspace{3mm}

{\it  Our aim is a generalization of the estimation (3), (4), (5) on the so - called Bilateral Grand Lebesgue Spaces $ BGL = BGL(\psi) = G(\psi), $ i.e. when } $ f(\cdot)
 \in G(\psi) \ $  {\it and to show the precision of obtained estimations by means
 of the constructions of suitable examples. } \par

 \vspace{3mm}

  We recall briefly the definition and needed properties of these spaces.
  More details see in the works \cite{Fiorenza1}, \cite{Fiorenza2}, \cite{Ivaniec1},
   \cite{Ivaniec2}, \cite{Ostrovsky1}, \cite{Ostrovsky2}, \cite{Kozatchenko1},
  \cite{Jawerth1}, \cite{Karadzov1} etc. More about rearrangement invariant spaces
  see in the monographs \cite{Bennet1}, \cite{Krein1}. \par

\vspace{3mm}

For $a$ and $b$ constants, $1 \le a < b \le \infty,$ let $\psi =
\psi(p),$ $p \in (a,b),$ be a continuous positive
function such that there exists a limits (finite or not)
$ \psi(a + 0)$ and $\psi(b-0),$  with conditions $ \inf_{p \in (a,b)} > 0 $ and
 $\min\{\psi(a+0), \psi(b-0)\}> 0.$  We will denote the set of all these functions
 as $ \Psi(a,b). $ \par

The Bilateral Grand Lebesgue Space (in notation BGLS) $  G(\psi; a,b) =
 G(\psi) $ is the space of all measurable
functions $ \ f: R^d \to R \ $ endowed with the norm

$$
||f||G(\psi) \stackrel{def}{=}\sup_{p \in (a,b)}
\left[ \frac{ |f|_p}{\psi(p)} \right], \eqno(6)
$$
if it is finite.\par
 In the article \cite{Ostrovsky2} there are many examples of these spaces.
 For instance, in the case when  $ 1 \le a < b < \infty, \beta, \gamma \ge 0 $ and

 $$
 \psi(p) = \psi(a,b; \beta, \gamma; p) = (p - a)^{-\beta} (b - p)^{-\gamma};
 $$
we will denote
the correspondent $ G(\psi) $ space by  $ G(a,b; \beta, \gamma);  $ it
is not trivial, non - reflexive, non - separable
etc.  In the case $ b = \infty $ we need to take $ \gamma < 0 $ and define

$$
\psi(p) = \psi(a,b; \beta, \gamma; p) = (p - a)^{-\beta}, p \in (a, h);
$$

$$
\psi(p) = \psi(a,b; \beta, \gamma; p) = p^{- \gamma} = p^{- |\gamma|}, \ p \ge h,
$$
where the value $ h $ is the unique  solution of a continuity equation

$$
(h - a)^{- \beta} = h^{ - \gamma }
$$
in the set  $ h \in (a, \infty). $ \par

 The  $ G(\psi) $ spaces over some measurable space $ (X, F, \mu) $
with condition $ \mu(X) = 1 $  (probabilistic case)
appeared in \cite{Kozatchenko1}.\par
 The BGLS spaces are rearrangement invariant spaces and moreover interpolation spaces
between the spaces $ L_1(R^d) $ and $ L_{\infty}(R^d) $ under real interpolation
method \cite{Carro1}, \cite{Jawerth1}. \par
It was proved also that in this case each $ G(\psi) $ space coincides
with the so - called {\it exponential Orlicz space,} up to norm equivalence. In others
quoted publications were investigated, for instance,
 their associate spaces, fundamental functions
$\phi(G(\psi; a,b);\delta),$ Fourier and {\it singular} integral operators,
conditions for convergence and compactness, reflexivity and
separability, martingales in these  spaces, etc.\par

{\bf Remark 1.} If we introduce the {\it discontinuous} function

$$
\psi_r(p) = 1, \ p = r; \psi_r(p) = \infty, \ p \ne r, \ p,r \in (a,b)
$$
and define formally  $ C/\infty = 0, \ C = \const \in R^1, $ then  the norm
in the space $ G(\psi_r) $ coincides with the $ L_r $ norm:

$$
||f||G(\psi_r) = |f|_r.
$$
Thus, the Bilateral Grand Lebesgue spaces are direct generalization of the
classical exponential Orlicz's spaces and Lebesgue spaces $ L_r. $ \par

 The function $ \psi(\cdot) $ may be generated as follows. Let $ \xi = \xi(x)$
be some measurable function: $ \xi: X \to R $ such that $ \exists  (a,b):
1 \le a < b \le \infty, \ \forall p \in (a,b) \ |\xi|_p < \infty. $ Then we can
choose

$$
\psi(p) = \psi_{\xi}(p) = |\xi|_p.
$$

 Analogously let $ \xi(t,\cdot) = \xi(t,x), t \in T, \ T $ is arbitrary set,
be some {\it family } $ F = \{ \xi(t, \cdot) \} $ of the measurable functions:
$ \forall t \in T  \ \xi(t,\cdot): X \to R $ such that
$$
 \exists  (a,b): 1 \le a < b \le \infty, \ \sup_{t \in T} \
|\xi(t, \cdot)|_p < \infty.
$$
Then we can choose

$$
\psi(p) = \psi_{F}(p) = \sup_{t \in T}|\xi(t,\cdot)|_p.
$$
The function $ \psi_F(p) $ may be called as a {\it natural function} for the family $ F. $ This method was used in the probability theory, more exactly, in
the theory of random fields, see \cite{Ostrovsky1}. \par

 The BGLS norm estimates, in particular, Orlicz norm estimates for
measurable functions, e.g., for random variables are used in the theory of
Partial Differential Equations \cite{Fiorenza1}, \cite{Ivaniec1}, theory of
probability in Banach spaces \cite{Ledoux1}, \cite{Kozatchenko1},
\cite{Ostrovsky1}, in the modern non-parametrical statistics, for
example, in the so-called regression problem \cite{Ostrovsky1}.\par

The article is organized as follows. In the second section we obtain
the main result: upper bounds for {\it integral} operators in the Bilateral
Grand Lebesgue spaces. In the next sections we investigate some one-dimensional
classical {\it singular } operators in BGLS: we obtain the upper estimations for
its norm and consider some examples in order to show the sharpness of upper
estimations.\par

 Fourth section contains the asymptotically exact BGLS norm calculation for the
 integral operators with homogeneous kernel. In the sixth sections we investigate
 some integral operators over the spaces with weight. \par

 The $ 7^{th} $ section contains the investigation of Fourier integral operators
in BGLS spaces with exact embedding constant calculations. In the eight section we
obtain the exact norm estimation for the celebrated Riesz operator.\par
 In the $ 9^{th} $ section we investigate the case when the measures $ \mu_{1,2} $
are pure atomic with unit measure of all points ("discrete case"). \par

 The last section contains some concluding remarks and generalizations. \par

\vspace{3mm}

 We use symbols $C(X,Y),$ $C(p,q;\psi),$ etc., to denote positive
constants along with parameters they depend on, or at least
dependence on which is essential in our study. To distinguish
between two different constants depending on the same parameters
we will additionally enumerate them, like $C_1(X,Y)$ and
$C_2(X,Y).$ The relation $ g(\cdot) \asymp h(\cdot), \ p \in (A,B), $
where $ g = g(p), \ h = h(p), \ g,h: (A,B) \to R_+, $
denotes as usually

$$
0< \inf_{p\in (A,B)} h(p)/g(p) \le \sup_{p \in(A,B)}h(p)/g(p)<\infty.
$$
The symbol $ \sim $ will denote usual equivalence in the limit
sense.\par
 We will denote as ordinary the indicator function
$$
I(x \in A) = 1, x \in A, \ I(x \in A) = 0, x \notin A;
$$
here $ A $ is a measurable set.\par
 All the passing to the limit in this article may be proved by means
 of Lebesgue dominated convergence theorem.\par

\bigskip

\section{Main result: norm estimations for regular integral operators.}

\vspace{3mm}

Let $ \psi(\cdot) \in \Psi(a,b), \ $ where $ 1 \le a < b.  $
 Let $ T $ be the operator described in the first section. We introduce the function

 $$
 \nu(p) = \inf_{q \in (a,b) }  \left\{ |T|_{p,q} \cdot \psi(q) \right\}. \eqno(7)
 $$
 Then the function $ \nu(\cdot) $ belong to the space $ G(\psi;c,d), $ where
 $$
 c = \inf \{p, p \ge 1, \nu(p) < \infty \}; \ d = \sup \{p, p \ge 1, \nu(p) < \infty \}.
 $$
  For the function $ \nu(\cdot) $ is true the elementary estimation

 $$
 \nu(p) \le \min_j \inf_{q \in (a,b) }  \left\{ |K|^{(j)}_{p,q} \cdot \psi(q) \right\}. \eqno(8)
 $$

\vspace{3mm}

{\bf Theorem 1.} Let $ f \in G(\psi), \ \psi \in \Psi(a,b) $  and
let $ T $ be the operator described in the first section.  Then

$$
||T \ f||G(\nu) \le || f ||G(\psi).  \eqno(9).
$$
 On the other words, the operator $ T $ is bounded from the space $ G(\psi) $ into the {\it another} space $ G(\nu). $ \par

\vspace{3mm}

{\bf Proof} of the theorem 1 is very simple. Denote for the simplicity
$ u = T f. $  We suppose $ f(\cdot) \in G(\psi); $ otherwise is nothing to prove.\par

 We can assume without loss of generality that $ ||f ||G(\psi) = 1; $
this means that

$$
\forall q \in (a,b) \ \Rightarrow |f|_q \le \psi(q).
$$

Using the definition of the norm $ |T|_{p,q} $ we obtain:

$$
|u|_p \le |T|_{p,q} \cdot  \psi(q), \ q \in (a,b).
$$

The assertion of theorem 1 follows after
the dividing over the  $ \nu(q), $ tacking the minimum  over
$ q, \ q \in (a,b) $ and on the basis of the definition of the $ G(\psi) $ spaces.\par
\hfill $\Box$ \\

{\bf Corollary 1.}  We consider separately the case when $ T \in UC; $  then the
assertion (9) of theorem 1 may be rewritten as follows:

$$
|T \ f|_{w(p)} \le \nu(p) \cdot |f|_p. \eqno(10)
$$

\vspace{3mm}

 We investigate in the remainder part of this section the {\it sufficient} conditions for {\it compactness } of the operator $ T $  from the one  Bilateral Grand Lebesgue Space  to the other  Bilateral Grand Lebesgue  Space.\par
 We assume here that both the measures $ \mu_1,\Sigma_1 $ and  $ \mu_2, \Sigma_2 $ are separable  relative the distances correspondingly

$$
\rho_1(A_1,A_2) = \mu_1(A_1 \setminus A_2) + \mu_1(A_2 \setminus A_1), A_1,A_2 \in \Sigma_1;
$$

$$
\rho_2(B_1,B_2) = \mu_2(B_1 \setminus B_2) + \mu_2(B_2 \setminus B_1), B_1,B_2 \in \Sigma_2.
$$
 It follows from the famous Kondrashov theorem \cite{Kantorovicz1}, chapter 9,
 section 3 that if $ K^{(3)}_{p,q} < \infty, $ then the operator $ T $ is compact as
 the  linear operator from the space $ L_{q,\mu_2} $ into the space  $ L_{p,\mu_1}. $ \par
 For example, it is true for the so-called embedding Sobolev's operator. \par
 But in the case of BGLS spaces still the embedding operator is bounded and is not
 compact  \cite{Ostrovsky2}. \par
{\bf Theorem 2.} Let $ \psi(\cdot) \in G(\psi;a,b), $ and let
$$
\max(\lim_{p \to c+0} \nu(p), \lim_{p \to d-0} \nu(p)) = \infty. \eqno(11)
$$
Let $  \zeta(\cdot) $ be arbitrary function from the space $ G(\psi; c,d) $ such that $ \zeta(\cdot) >> \nu(\cdot). $ Then the operator $ T $ is compact operator from the space $  G(\psi) $  into the space $ G(\zeta). $ \par
 {\bf Proof. } It follows from the  Kondrashov's theorem that the operator $ T $ is compact operator from the arbitrary space $ L_{p, \mu_2}, \ p \in (a,b) $ into the arbitrary space $ L_{q, \mu_1}, \ q \in (c,d). $ The proposition of theorem 2 follows from the main result  of an article \cite{Ostrovsky2}.\par
\hfill $\Box$ \\
\bigskip

\section{Classical one-dimensional singular operators}

\vspace{3mm}
 We consider in this section some classical one-dimensional integral operators
with homogeneous kernel of degree $ = - 1 $ of Hardy-Littlewood-Young type.
We intend to obtain the {\it exact } values of embedding constants or as a minimum to obtain the low bounds for the norms of the operator $ T $ in the BGLS spaces.\par

In detail, we consider the one-dimensional singular integral operators in the space
$ X = (0,\infty) $ of a view

$$
T_H f(x) = \int_0^{\infty} K(x,y) \ f(y) \ dy, \eqno(12)
$$
where the kernel $ K(\cdot,\cdot) $ is presumed to be non-negative and homogeneous of degree -1; this means that
$$
K(x,y) = x^{-1} \ H(y/x), \ x,y > 0,
$$
where $ H(\cdot) $ is some measurable function. \par

We define for the values $ p \in (1,\infty) $  the function

$$
\phi(p) = \phi_H(p) \stackrel{def}{=} \int_0^{\infty} z^{ - 1/p} H(z) \ dz \eqno(13)
$$
and suppose in this subsection the {\it finiteness} of the function $ \phi(p) $ for
all the values s$ p \in (1,\infty). $ \par
 The famous result belonging to Hardy and Littlewood \cite{Hardy1}, p. 93-114
 states that

 $$
 |T_H|_{p,p} \le \phi(p).
 $$
and the last inequality is exact under some additional conditions. We can therefore
use theorem 1. Namely, we define the arbitrary  function $ \psi(\cdot) $ from the set
$ G\Psi(1,\infty) $ the new function
$$
\psi_{(\phi)}(p) = \phi(p) \psi(p).
$$
 It follows from theorem 1 that

$$
||T_H \ f||G(\psi_{(\phi)}) \le V_H \ ||f||G(\psi), \ V_H \le 1. \eqno(14)
$$
We will prove that under some additional conditions the exact value for the
 constant $ V_H $ in the inequality (14) is equal to one.\par
 Before the formulating the main result of this section, we need to establish
 some preliminary lemmas. \par

{\bf Lemma 1a.} Assume that the (measurable non-negative) function $ H(\cdot) $
satisfies the following conditions:
$$
\exists \delta_0 > 0, \ \int_1^{\infty} z^{-1  +\delta_0}
H(z) dz < \infty;  \eqno(15a)
$$
$$
\exists H_- \in (0,\infty), \ z \to 0+ \Rightarrow  H(z) =
$$
$$
H_- \ |\log z|^{\beta_1} \ S_1(|\log z|) ( 1 + 0(|\log z|^{-\gamma}|)), \eqno(15b)
$$
where  $ H_- = \const \in (0,\infty), \beta_1 = \const \ge 0;  \
S_1(z) $  is  non-negative  continuous slowly varying as $ z \to \infty $ functions:
$$
\forall l > 0 \ \Rightarrow \lim_{z \to \infty} S_1(l z)/S_1(z) = 1.
$$
 We propose that as $ p \to 1+0 $
$$
\phi_H(p) \sim H_-  (p-1)^{-\beta_1 +1} \ S_1(1/(p-1)) \ \Gamma(\beta_1+1). \eqno(16)
$$

{\bf Lemma 1b.} Assume that the (measurable non-negative) function $ H(\cdot) $
satisfies the following conditions:
$$
\exists \delta_0 > 0, \ \int_0^1 z^{-\delta_0} H(z) dz < \infty;  \eqno(17a)
$$
$$
\exists H_+ \in (0,\infty), \ z \to 0+ \Rightarrow  H(z) =
$$
$$
H_+ \ z^{-1} |\log z|^{\beta_2} \ S_2(|\log z|) ( 1 + 0(|\log z|^{-\gamma}|)), \eqno(17b)
$$
where  $ H_+ = \const \in (0,\infty), \beta_2 = \const \ge 0;  \
S_2(z) $  is  non-negative  continuous slowly varying as $ z \to \infty $ functions:
$$
\forall l > 0 \ \Rightarrow \lim_{z \to \infty} S_2(l z)/S_2(z) = 1.
$$
 We propose that as $ p \to \infty $
$$
\phi_H(p) \sim H_+ \ p^{\beta_2 +1} \ S_2(p) \ \Gamma(\beta_2+1). \eqno(18)
$$
{\bf Proofs}.  Proof of lemma 1a. We write first of all the partition
$$
\phi_H(p) = \int_0^1 z^{-1/p} H(z) dz + \int_1^{\infty} z^{-1/p} H(z) dz  =
I_1 + I_2,
$$
and note that the second integral is uniformly bounded for the values $ p $ nearest
to $ 1+0 $ by virtue of the condition (15a). \par
The assertion of lemma 1a  it follows from the asymptotical equality: as
$ p \to 1 + 0 $

$$
\int_0^1 z^{-1/p } \ |\log z|^{\beta_1} \ S_1(\log z) \ dz \sim
$$

$$
(p-1)^{-\beta_1 -1} \  S \left( \frac{1}{p-1} \right) \
\Gamma(\beta_1+1), \ \beta_1 \ge 0.
$$
The proof of the lemma 1.b may be obtained analogously. \par
\vspace{3mm}

{\bf Theorem 3.} Suppose the function $ H(\cdot) $ satisfies either the conditions
 of lemma 1a or the conditions of lemma 1b. Then
$$
||T_H \ f||G(\psi_{(\phi)})  \le V_H \ ||f||G(\psi), \eqno(19)
$$
where the exact value of constant $ V_H $ is equal to one.\par
{\bf Proof } of the theorem 3. \par

{\bf A.} Let us consider and investigate the following value:
$$
Z \stackrel{def}{=} \sup_{\psi \in G\Psi(1,\infty)}
||T_H||(G(\psi) \to G(\psi_{(\phi)})=
$$
$$
\sup_{\psi \in G\Psi(1,\infty)} \sup_{f \in G\psi}
\frac{ ||T_H [f]||G(\psi_{(\phi)})}{||f||G(\psi)} =
$$

$$
\sup_{\psi \in G\Psi(1,\infty)} \sup_{f \in G\psi}
\frac{\sup_p [|T_H[f]|_p/(\psi(p) \phi(p))]}{\sup_p [|f|_p/\psi(p)]}. \eqno(19)
$$

 Note that if we choose the function $ \psi(\cdot) $ as follows:
 $$
 \psi(p) = |f|_p,
 $$
i.e. the {\it natural } choice of a function $ \psi(\cdot), $ obviously if the
function $ f $ belong to the set $ L(1,\infty), $ we  obtain the following {\it low
estimation } for the value $ Z: $

$$
Z \ge \sup_{f \in L(1,\infty)}\sup_{p \in (1,\infty) }
\frac{|T_H f|_p}{\phi(p) \ |f|_p }. \eqno(20)
$$

{\bf B.} It remains only to prove the low bound for the value $ V_H. $ Let us consider at first the case when the function $ H(\cdot) $ satisfies the conditions of lemma 1a.\par
Let us choose the function, more exactly, the family of a functions

$$
g_{\Delta}(x) = |\log x|^{\Delta} I(x \in (0,1)),
$$
here $ \Delta = \const \ge 2 $ is some fixed number. \par
 We have:

$$
|g_{\Delta}|_p = \sqrt[p]{\Gamma(\Delta p + 1)}.
$$
 Further, we denote $ v_{\Delta} = T_H g_{\Delta}. $ We obtain as $x \to 0:$

 $$
 v_{\Delta}(x) = \int_0^1 x^{-1} \ |\log y|^{\Delta} \ H(y/x) \ dy \sim
 $$

 $$
\frac{ H_+}{\Delta+\beta_2 + 1} \ |\log x|^{\beta_2+\Delta+1}S_2(|\log x|).
 $$
 We obtain after some calculation

 $$
 \overline{\lim}_{p \to \infty} \frac{|v_{\Delta}|_p }{\phi(p) \ |g_{\Delta}|_p } \ge
\exp(-\beta_2 - 1) \ \left( 1+ \frac{\beta_2 + 1 }{\Delta} \right)^{\Delta}.
$$
Therefore,

$$
Z \ge \lim_{\Delta \to \infty}
\exp(-\beta_2 - 1) \ \left( 1+ \frac{\beta_2 + 1 }{\Delta} \right)^{\Delta}=1.
$$
This completes the proof of theorem 3 in the case when the function $ H(\cdot) $
satisfies the conditions of lemma 1b. \par
{\bf C.} The second case. i.e.when the function $ H(\cdot) $ satisfies the conditions
of lemma 1a, the proof provided analogously by the consideration  the family of a functions

$$
f_{\Delta}(x) = x^{-1} \ (\log x)^{\Delta}  \ I(x > 1).
$$

\vspace{3mm}

 We consider further in this section the integral operators with homogeneous kernel
 of degree -1 which does not satisfy the conditions of theorem 3. For instance, let
 the operator $ T_{r,s}[f] $ has a view:

 $$
 T_{r,s}[f](x) = \int_0^x \frac{(x-y)^{r-1} f(y) \ dy  }{x^s y^{r-s}},
 $$
 where $ r,s = \const, \ 0 < r < s+1. $ Here
 $$
  H(z)  = z^{-(r-s)} \ (1-z)^{r-1} \ I(z \in (0,1))
 $$
and correspondingly {\it only for the values} $ p > 1/(1+ s - r)  $

 $$
 \phi_{r,s}(p) = B( 1 - 1/p -(r-s), r ) =
\frac{\Gamma(1-1/p - (r-s)) \ \Gamma(r)}{\Gamma(1-1/p + s)},
$$
where $ B(\cdot,\cdot) $ denotes usually Beta-function. \par
 Recall that we consider here the case when $ p > 1; $ so we assume that
 $$
 p > \max(1, 1/(1+s-r)),
 $$
in the contradiction to the theorem 3.\par
 The family of the operators  $ \{ T_{r,s}[f](x) \} $ contains the classical Rieman's 
fractional integral operator, up to multiplicative constant and the power factor $ x^t,$
which does not dependent on the variable $ y: $
$$
R_{(r)}[f](x) \ \frac{1}{\Gamma(r)} \int_0^x (y-x)^{r-1} \ f(y) \ dy.
$$ 
 Analogously may be considered the fractional integral in the Weil's sense: 

$$
W^{(r)}[f](x) \ \frac{1}{\Gamma(r)} \int_x^{\infty} (y-x)^{r-1} \ f(y) \ dy.
$$

\vspace{3mm}

So, we will consider the operators $ T_H[f] $ with homogeneous kernel
$$
 K(x,y) = x^{-1}H(y/x), \ \degree( K) = - 1
$$
where the function $ H = H(z) $ does not satisfy the conditions of lemmas 1a or
lemma 1b. Indeed, we consider as a functions $ H(\dot) $ the functions of a view

$$
{\bf A} \ H(z) \sim H_- \cdot z^{\alpha - 1} |\log z|^{\beta_1} S_1(z), \ z \to 0+
$$
or

$$
{\bf B} \ H(z) \sim H_+ \cdot z^{\alpha - 1} |\log z|^{\beta_1} S_2(z), \ z \to \infty.
$$
In both the cases $ \alpha = \const \in (0,1) $  and as before $ S_1(\cdot), S_2(\cdot) $
are non-negative continuous  slowly varying as $ z \to \infty $ functions. \par
 Note that in considered cases {\bf A } and {\bf B} the function $ \phi(p) $ does not
 exists in the whole axis $ (1,\infty), $ in contradiction to the cases considered in
 the lemmas 1a and 1b. Namely, in the first case {\bf A}

 $$
 \phi_H(p) < \infty \leftrightarrow p > 1/\alpha.
 $$
and in the case of the condition {\bf B}

 $$
 \phi_H(p) < \infty \leftrightarrow p \in (1, 1/\alpha).
 $$

{\bf Lemma 2a.}

 Assume that the (measurable non-negative) function $ H(\cdot) $
satisfies the following conditions:
$$
\exists \delta_0 > 0, \ \int_1^{\infty} z^{-1  +\delta_0}
H(z) dz < \infty;
$$
$$
\exists H_- \in (0,\infty), \ z \to 0+ \Rightarrow  H(z) =
$$
$$
H_- \ z^{\alpha-1} \ |\log z|^{\beta_1} \ S_1(|\log z|) ( 1 + 0(|\log z|^{-\gamma}|)),
$$
where  $ H_- = \const \in (0,\infty), \beta_1 = \const \ge 0;  \
S_1(z) $  is  non-negative  continuous slowly varying as $ z \to \infty $ functions.

 We propose that as $ p \to 1/\alpha+0 $
$$
\phi_H(p) \sim H_- \ (p-1/\alpha)^{-\beta_1 +1} \ S_1(1/(p-1/\alpha)) \ \Gamma(\beta_1+1). $$

{\bf Lemma 2b.}
 Assume that the (measurable non-negative) function $ H(\cdot) $
satisfies the following conditions:
$$
\exists \delta_0 > 0, \ \int_0^1 z^{-\delta_0} H(z) dz < \infty;
$$
$$
\exists H_+ \in (0,\infty), \ z \to 0+ \Rightarrow  H(z) =
$$
$$
H_+ \ z^{-1 - \alpha} |\log z|^{\beta_2} \ S_2(|\log z|) ( 1 + 0(|\log z|^{-\gamma}|)),
$$
where  $ H_+ = \const \in (0,\infty), \beta_2 = \const \ge 0;  \
S_2(z) $  is  non-negative  continuous slowly varying as $ z \to \infty $ function:

 We propose that as $ p \to 1/\alpha - 0 $
$$
\phi_H(p) \sim H_+ \ (1/\alpha-p)^{-\beta_2 +1}
 \ S_2(1/\alpha - p) \ \Gamma(\beta_2+1).
$$
 {\bf Proofs} are analogous to the lemma 1a and may be omitted.\par

 \vspace{3mm}

{\bf Theorem 4.} Suppose the function $ H(\cdot) $ satisfies either the conditions
 of lemma 2a or the conditions of lemma 2b. Then
$$
||T_H \ f||G(\psi_{(\phi)})  \le V_H \ ||f||G(\psi), \eqno(20)
$$
where the exact value of constant $ V_H $ is equal to one.\par
{\bf Proof } of the theorem 4 is at the same as the proof of theorem 3. \par

\hfill $\Box$ \\

\bigskip

\section{ Singular operators with arbitrary degree of kernel, with
generalization}

\vspace{3mm}

{\bf A.} We consider in this subsection the {\it family } of operators of a view

$$
T_{\lambda}f(x) = \int_0^{\infty} \frac{ f(y) \ dy  }{ (x+y)^{\lambda}}.
$$

Here $ \lambda = \const \in (0,1), \ X = (0,\infty) $  and we denote
$ u(x) = T_{\lambda}f(x). $ \par

 It is known \cite{Mitrinovich1}, p. 215 that

$$
|T_{\lambda}f|_{p/\lambda} \le
\left[ \frac{ \pi }{\sin(\pi/p) }  \right]^{\lambda} |f|_{q/(q-\lambda)},
$$
where $ p \in (\lambda, \infty) $ and as ordinary $ q = p/(p-1). $ \par
 The last relation may be rewritten as follows.

$$
|T_{\lambda}f|_p \le \left[ \frac{ \pi }{\sin(\pi/(p \lambda)) }  \right]^{\lambda} \cdot |f|_{p/(p(1-\lambda)+1) },
$$
but for the values $ p $ from the interval $ p \in (1/\lambda, \infty). $ \par
 In order to formulate the main result of this section, we need to introduce some functions. Let $ h =  h(\lambda), \ \lambda \in (0,1) $ be well-known
entropy  function

$$
h(\lambda) = -\lambda \log(\lambda) \ - \ (1-\lambda) \log(1-\lambda).
$$
Note that the function $ h(\cdot) $ may be continued on the whole closed interval $ [0,1] $ as the continuous function, as long as

$$
\lim_{\lambda \to 0+} h(\lambda) = \lim_{\lambda \to 1-0} h(\lambda) = 0.
$$
 Further, we define the function:

$$
L(\lambda): = \max \left\{ \lambda \Gamma^{\lambda}(1+1/\lambda),
e^{h(\lambda)} \right\},
$$

$$
V(\lambda,p) = \sup_{f \in L(1/\lambda, \infty), f \ne 0}
\left\{{\frac{|T_{\lambda}f|_{p/\lambda}}{|f|_{p/( p(1-\lambda) + \lambda) }}}
 \cdot \left[ \frac{\sin(\pi/p)}{\pi} \right]^{\lambda} \right\},
$$
$ p \in (1,\infty); $

$$
\overline{V}(\lambda)= \sup_{p \in (1,\infty)} V(\lambda,p).
$$
Note that
$$
 W \stackrel{def}{=} \sup_{\lambda \in (0,1)} \overline{V}(\lambda) = 1.
$$

 Let also $ \psi(\cdot) $ be arbitrary function from the space
$ G\Psi(1, 1/(1-\lambda)). $ We define for any value $ \lambda \in (0,1) $

$$
\psi_{(\lambda)}(p):=  \left[ \frac{\pi}{\sin(\pi/(p \lambda))} \right]^{\lambda}
\cdot \psi \left( \frac{ p }{ p(1-\lambda) + 1 } \right).
$$

{\bf Theorem 5.}

$$
||T^{(1)} \ f||G(\psi^{(1)}) \le V_1 \ ||f||G(\psi), \eqno(21)
$$
where the exact value of constant $ V_1 =V_1(\lambda) $  lies in the closed interval
$$
V_1 \in [L(\lambda), 1].
$$
{\bf Proof.} The upper bound for $ V_1 $ it follows from the $ L_p \to L_p $ estimations
for considered operator on the basis of theorem 1. To obtain the low bound
$  L(\lambda), $ we construct two examples. \par
{\sc First example.} Let us consider the function

$$
f(x) = x^{-1} I(x>1),
$$
then

$$
|f|_Q =(Q-1)^{-1/Q}, \ Q > 1.
$$
and as $ x \to \infty $
$$
u_{\lambda}(x):= \int_1^{\infty} \frac{dy}{y(x+y)^{\lambda}} = x^{-\lambda}
\int_{1/x}^{\infty} \frac{dz}{z(1+z)^{\lambda}} \sim
$$

$$
x^{-\lambda} \int_{1/x}^1 z^{-1} dz = x^{-\lambda} \log x;
$$
and we compute ad $ p \to 1+0 $
$$
|u_{\lambda}|_p \sim \frac{ \Gamma^{1/p}(p+1)}{(p\lambda - 1)^{1+1/p}}.
$$
The first example give us the following low estimation for  the value $ V_1: $
$$
V_1 \ge  \lambda \Gamma^{\lambda}(1+1/\lambda).
$$

{\sc Second example.} We consider here the function

$$
g(x) = x^{-(1-\lambda)} I(x \in (0,1)).
$$
We get as $ p \to \infty:$
$$
|g|_Q  \sim \left( \frac{1-\lambda}{\lambda} \right)^{1-\lambda}p^{1-\lambda};
$$

$$
v(x):= \int_0^1 \frac{y^{-(1-\lambda)} \ dy}{(x+y)^{\lambda}} \sim |\log x|,
$$
$$
|v|_p^p \sim \int_0^1 |\log x|^p dx = \Gamma(p+1).
$$

The second example give us the next low estimation for  the value $ V_1: $
$$
V_1 \ge e^{h(\lambda)}.
$$
Thus,
$$
V_1 \ge \max \left(  e^{h(\lambda)}, \Gamma^{\lambda}(1+1/\lambda) \right),
 \ \lambda \in (0,1). \eqno(22)
$$
\hfill $\Box$ \\

\vspace{3mm}

{\bf B.} We consider here  the {\it family } of operators of a view

 $$
u(x) = I_{\alpha,\beta,\lambda}[f](x) = I[f](x) = |x|^{-\beta}
\int_{R^d} \frac{f(y) \ |y|^{-\alpha} \ dy}{|x - y|^{\lambda}}; \eqno(23)
 $$
here $ \alpha,\beta \ge 0, \alpha + \lambda < d.  $ \par
 We define the function $ q=q(p) $ as follows:
$$
1+\frac{1}{q} = \frac{1}{p} + \frac{\alpha+\beta+\gamma}{d}.
$$
We will denote the set of all such a values $ (p,q) $  as $ G(\alpha,\beta,\lambda) $ or
for simplicity $ G = G(\alpha,\beta,\lambda). $  \par
{\sc Further we will suppose in this subsection that }
$ (p,q) \in G(\alpha,\beta,\lambda) = G. $ \par

  We denote also
$$
p_-:=\frac{d}{d-\alpha}, \ p_+:= \frac{d}{d-\alpha-\lambda};
$$
and correspondingly
$$
q_-:= \frac{d}{\beta+\lambda}, \ q_+ := \frac{d}{\beta},
$$
where in the case $ \beta=0 \ \Rightarrow q_+:= + \infty; $
$$
\kappa = \kappa(\alpha,\beta,\lambda) := (\alpha + \beta + \lambda)/d.
$$
 Let $ \psi(\cdot) \in G\Psi(p_-,p_+); $ we define

 $$
 \psi_{\alpha,\beta,\lambda} = \left[(p-p_-)(p_+ - p) \right]^{-\kappa} \ \psi(p).
 \eqno(24)
 $$

{\bf Theorem 6.} There exist a two positive finite constants
$ C_1(d; \alpha,\beta, \lambda), \ C_2(d; \alpha,\beta, \lambda) $ for which
the minimal value of a constant $ \le V_{\alpha,\beta,\lambda} $ from the inequality
$$
||I_{\alpha,\beta,\lambda}[f]||G \left(\psi_{\alpha,\beta,\lambda}\right) \le V_{\alpha,\beta,\lambda} \ ||f||G(\psi)
$$
lies in the closed interval

$$
 V_{\alpha,\beta,\lambda} \in [C_1(d; \alpha,\beta, \lambda),
C_1(d; \alpha,\beta, \lambda)]. \eqno(25)
$$
{\bf Proof} follows immediately from the main results of papers \cite{Ostrovsky3}, \cite{Ostrovsky4}. There are described, for instance, correspondent examples.\par
 Note that in the case $ \alpha = \beta = 0 $ we can  conclude that
 $$
\frac{C_3(d)}{d-\lambda} \le V_{0,0,\lambda} \le \frac{C_4(d)}{ d - \lambda},
 $$
see \cite{Ostrovsky3}. \par
\hfill $\Box$ \\
 \bigskip

\section{ Singular operators of Hardy type and its weight generalizations}

\vspace{3mm}
We consider in this section the operators of a view

$$
U_{\lambda}[f](t) = t^{\lambda-1} \int_0^t s^{-\lambda} \ f(s) \ ds.
$$
Here $ X = R_+, \ \lambda = \const \in (-\infty,1), $ and we equip all the Borelian
sets of the space $ X $ by the measure
$$
\nu(A) = \int_A \frac{dt}{t}.
$$
 This operators play a very important role in the theory of operators interpolation
 \cite{Bennet1}, chapter 3, section 1. \par
  There it is proved the following $ L_{p,\nu} $ estimation:

$$
|U_{\lambda}[f]|_{p,\nu} \le \frac{1}{1-\lambda} \ |f|_{p,\nu},  \ p \in [1,\infty].
$$
Therefore, if we consider arbitrary function $ \psi \in G\Psi(1,\infty), $
we conclude on the basis of theorem 1 that

$$
||U_{\lambda}[f]||G(\psi) \le \frac{1}{1-\lambda} \ ||f||G(\psi). \eqno(26)
$$

{\bf Theorem 7}. The constant $ 1/(1-\lambda) $ in the inequality (26) is exact. \par
{\bf Proof } is at the same as before. Namely, let us consider the function
$$
f_0(t) = |\log t|^{-\Delta} \ I(t \in (0,1)), \ \Delta = \const \in (0,1);
$$
then
$$
|f_0|_{p,\nu} = (p \Delta - 1)^{-1/p}, \ p \in (1/\Delta, \infty),
$$
and we get as $ t \to 0+: $
$$
U_{\lambda}[f_0](t) = t^{\lambda -1}\int_0^t \ s^{-\lambda} \ |\log s|^{-\Delta} \ ds \sim
$$

$$
\frac{1}{1-\lambda} \ |\log t|^{-\Delta} = \frac{1}{1-\lambda} \ f(t).
$$
Therefore, we have as $ p \to 1/\Delta + 0 $
$$
|U_{\lambda}[f_0]|_p \sim \frac{1}{1-\lambda} |f_0|_p.
$$
This completes the proof of this theorem.\par
 Another version of Hardy's inequality see in \cite{Reyna1}, pp. 134-135:

$$
\left( \int_0^{\infty} t^{-p/p_1} \left(\int_0^t |f(s)| \ s^{1/p_1 - 1} \ ds \right) \ dt         \right)^{1/p} \le
$$

$$
\left(1/p_1 - 1/p \right)^{-1} \ \left( \int_0^{\infty} |f(s)|^p \ ds \right)^{1/p};
\eqno(27a)
$$

$$
\left( \int_0^{\infty} t^{-p/p_0} \left(\int_t^{\infty} |f(s)| \ s^{1/p_0 - 1} \ ds \right) \ dt \right)^{1/p} \le
$$

$$
\left(1/p - 1/p_0 \right)^{-1} \ \left( \int_0^{\infty} |f(s)|^p \ ds \right)^{1/p};
\eqno(27b)
$$
here $ 1 \le p_1 < p < p_0 < \infty. $ \par
 The inequalities (27a) and (27b) may be rewritten as follows.
$$
U_{\alpha}[f]|_p \le (\alpha - 1/p)^{-\alpha} |f|_p, \ p \in (1/\alpha, \infty),
$$

$$
U_{\alpha}[f]|_p \le (1/p - \alpha )^{-\alpha} |f|_p, \ p \in (1, 1/\alpha),
$$
but here $ \alpha \in (0,1). $ \par
 If the function $ \psi = \psi(p) $ belongs to the set $ G\Psi(1/\alpha,\infty), $
we introduce the auxiliary function

$$
\psi_{\alpha}(p) = (\alpha - 1/p)^{-\alpha} \ \psi(p).
$$
It follows from the theorem 1 that

$$
||U_{\alpha}[f]||G(\psi_{\alpha}) \le 1 \cdot ||f||G(\psi). \eqno(28)
$$

{\bf Theorem 8.} The constant $ 1 $ in the inequality (28) is exact. \par
{\bf Proof }  is at the same as in the theorem 8. Namely, we consider  a function
$$
f_{\Delta}(x) = x^{-\alpha} \ (\log x)^{\Delta} \ I(x > 1), \ \Delta = \const \ge 1;
$$
then we have as $ p \to 1/\alpha + 0 $ and $ x  \to \infty, \ x \ge 1: $
$$
|f_0|_p = (\alpha p - 1)^{-\Delta -1/p} \ \Gamma^{1/p}(\Delta p + 1)
 \sim (\alpha p - 1)^{-\Delta -\alpha} \ \Gamma^{\alpha}(\Delta/\alpha + 1);
$$

$$
u_{\Delta}(x) := x^{-\alpha} \int_0^x s^{\alpha-1} s^{-\alpha} (\log s)^{\Delta} ds =
$$

$$
(\Delta+1)^{-1} \ x^{-\alpha} \ (\log x)^{\Delta+1},
$$

$$
|u_{\Delta}(\cdot) |_p \sim \frac{1}{\Delta+1}
\frac{\Gamma^{\alpha}((\Delta+1)/\alpha +1 )}{(\alpha p -1)^{(\Delta +1)/\alpha}}.
$$
 Note that as $ p \to 1/\alpha + 0 $

$$
\overline{\lim}_{p \to 1/\alpha + 0}
\frac{|u_{\Delta}|_p}{|f_{\Delta}|_p} : (\alpha - 1/p)  \ge
e^{-1} \left( \frac{\Delta +1}{\Delta}  \right)^{\Delta}. \eqno(29)
$$
 The expression (29) tends to one as $ \Delta \to \infty. $ \par
This completes the proof of this theorem.\par
The case when $ p \in (1, 1/\alpha)  $ may be considered analogously, by mean of
an example

$$
g_{\Delta}(x) = x^{-\alpha} \ |\log x|^{\Delta} \ I(x \in (0,1)).
$$

\hfill $\Box$ \\

\bigskip

\section{ Singular operators over the spaces with weight}

\vspace{3mm}
 Let $ r $ be arbitrary constant positive number, $ r \ne 1. $
 Let us consider in this section integral operators of a view

$$
f_{(1)}[f](x) = \int_0^x \ f(t) \ dt, \ X = R_+, \ r > 1,
$$
or
$$
f^{(1)}[f](x) = \int_x^{\infty} f(t) dt, \ X = R_+, \ r < 1.
$$
i.e. we consider the integral of $ f(\cdot). $  \par
We equip the Borelian sets of the space $ X = R_+ $ by means of the  weight measure
$ \mu_r(A): $

$$
\mu_r(A) = \int_A x^{-r} \ dx
$$
and will denote for simplicity  the $ L_p $ norm of the function $ h: R_+ \to R $
in this space as $ |h|_{p,r}: $

$$
|h|_{p,r}:= \left[\int_0^{\infty} x^{-r} \ |h(x)|^p \ dx \right]^{1/p}.
$$
 In the book \cite{Hardy1}, p. 227 - 231 it is proved:

 $$
 \int_0^{\infty}x^{-r} |f^{(1)}(x)|^p dx \le
 \left( \frac{p}{|r-1|} \right)^p \cdot  \int_0^{\infty} x^{-r} (x f(x))^p dx.
 \eqno(30)
 $$

The inequality (30) may be rewritten as follows:

$$
|F|_{p,r} \le \frac{p}{|r-1|} \cdot |g|_{p,r},
$$
where

$$
F[g](x) = \int_0^x y^{-1} \ g(y) \ dy, \ r > 1;
$$
the case when  $ r < 1 $ and

$$
F[g](x) = \int_x^{\infty} y^{-1} \ g(y) \ dy
$$
is considered analogously  and will be omitted. \par
Therefore, in the considered case $ r < 1. $ \par
 Let $ \psi(\cdot) $ be any function from the set $ G\Psi(1,\infty); $ we define

 $$
 \psi_1(p) = p \ \psi(p).
 $$

{\bf Theorem 9.}

$$
||F[g]||G(\psi_1) \le (|r-1|)^{-1} \ ||g||G(\psi),
$$
where the value of the constant $ (|r-1|)^{-1} $ is the best possible. \par
{\bf Proof} of the low estimate is ordinary. Let us consider the example of a
function:

$$
g_{\Delta}(x) =  |\log x|^{\Delta} I(x \in (0,1)), \ \Delta = \const > 1;
$$
then
$$
|g_{\Delta}|^p_{p,r} = \int_0^1 x^{-r} |\log x|^{\Delta p} dx =
\frac{\Gamma(\Delta p + 1)}{|1-r|^{p \Delta + 1}},
$$
and for the values $ x $ from the interval $ x \in (0,1): $

$$
F[g_{\Delta}](x) = \int_0^x y^{-1} \ |\log y|^{\Delta} \ dy = (\Delta+1)^{-1} \
|\log x|^{\Delta+1};
$$

 $$
 |F[g_{\Delta}]|_{p,r} \ge (\Delta + 1)^{-1}
  \frac{\Gamma^{1/p}(( \Delta+1)p + 1) }{|1-r|^{\Delta+1+1/p }}.
 $$
 It is easy to calculate using Stirling's formula that as $ p \to \infty $

 $$
 \overline{\lim} \left[ |r-1| \ |F[g_{\Delta}]|_p /|g_{\Delta}|_p \right] \ge
 e^{-1} (1+1/\Delta)^{\Delta}.
 $$
 The assertion of our theorem follows from the equality

 $$
 \lim_{\Delta \to \infty} e^{-1} (1+1/\Delta)^{\Delta}= 1.
 $$

 \hfill $\Box$ \\

\bigskip

\section{ Fourier integral operators }

\vspace{3mm}

We consider here the classical Fourier integral operator of a view

$$
F[f](x) = F(x) = \int_{-\infty}^{\infty} \exp(itx) f(t) dt.
$$

It is known \cite{Titchmarsh1}, pp. 96-98, that
$$
\int_R |F(x)|^p \ |x|^{p-2} dx \le p^2/(p-1) \ \int_R |f(t)|^p \ dt, \ p \in (1,\infty).
\eqno(31)
$$
 The inequality (31) may be rewritten in the language of BGLS as follows. Let $ \psi $
 be arbitrary function from the class $ G\Psi(1,\infty). $ We define the new measure
 $ \mu(dx) = dx/x^2 $ and  introduce the new function
 $$
 \psi_F(p) = p^2 \ \psi(p)/(p-1),
 $$
then it follows from the theorem 1 that
$$
||F[f]||G(\psi_F, \mu) \le 1 \cdot ||f||G(\psi). \eqno(32)
$$

{\bf Theorem 10.} The constant 1 in the inequality (32) is the best possible.\par
{\bf Proof.} Let us consider the example function

$$
f_0(x) = x^{-1} \ I(x > 1),
$$
then
$$
|f_0|_p = (p-1)^{1/p}.
$$
Further, we have as $ t\to 0, \ t \in (0,1): $

$$
F[f_0](t) = \int_1^{\infty} \exp(itx) dx/x = \int_t^{\infty} \exp(iy) dy/y = C +
$$

$$
\int_t^1 \exp(iy) dy/y \sim |\log t|.
$$

Therefore, as $ p \to 1+0 $

$$
|F[f_0]|^p_{p,\mu} \sim \int_0^1 |\log t|^p t^{p-2} dt =
\int_0^{\infty} e^{-y(p-1)  } y^p dy = (p-1)^{-p-1} \ \Gamma(p+1);
$$

$$
|F[f_0]|_{p,\mu} \sim (p-1)^{-1-1/p}
$$
and we have as before choosing $ \psi(p) = |f_0|_p $

$$
 \sup_{p \in (1,\infty)}
\frac{|F[f_0]|_{p,\mu}}{|f_0|_p} \times \frac{p-1}{p^2} \ge
 \overline{\lim}_{p \to 1}
\frac{|F[f_0]|_{p,\mu}}{|f_0|_p} \times \frac{p-1}{p^2} = 1.
$$

\hfill $\Box$ \\

\bigskip

\section{Riesz singular integral operator}

\vspace{3mm}

The operator of a view

$$
R_j[f](x) = c(n) \ p.v. \int_{R^n} \frac{(x_j-y_j) f(y) \ dy}{|x-y|^{n+1}}
$$
is named Riesz transform, or Riesz operator. Here $ n=1,2,\ldots;
j=1,2,\ldots,n;  $

$$
c(n) = \pi^{-(n+1)/2} \Gamma((n+1)/2).
$$
Note that the general case of singular integral operators acting in the BGLS spaces
without exact constant calculation is considered in \cite{Ostrovsky2}.\par
  In the one-dimensional case, i.e. when $ n=1 $ the Riesz transform coincides with the so-called Hilbert transform
$$
Rf(x) = Hf(x)\stackrel{def}{=} \pi^{-1} \ p.v. \int_{-\infty}^{\infty} \frac{f(y) \ dy }
{x-y}.
$$
 The exact value of the $ L_p \to L_p $ norm of the Riesz operator does not dependent on the dimension $ n. $   Namely, we have for the values $ p $ from
the {\it open} interval $ p \in (1,\infty): $

$$
h(p) \stackrel{def}{=} ||R_j||(L_p \to L_p) = \tan(\pi/(2p)), \ p \in (1,2]
$$
and

$$
h(p) = 1/\tan(\pi/(2p)), \ p \in [2,\infty),
$$
see \cite{Pichorides1}, \cite{Iwaniec3},  chapter 12, section 12.1.  \par
Note that as $ p \to 1 + 0 \ h(p) \sim 1/(p-1) $ and as $ p \to \infty \ h(p) \sim p.$ \par
 At the same estimation for the $ L_p \to L_p $ norm is true for the Hilbert transform
 in the case $ X = (-\pi, \pi),  $  where the operator $ H $ is defined as follows:

 $$
 H_{\pi}[f](x) = (2 \pi)^{-1} \ p.v. \ \int_{-\pi}^{\pi}
 \frac{f(x-y)}{\tan(y/2)} dy, \ x-y = x-y \mod(2 \pi).
 $$
 Note also that if
 $$
 f(x) = \frac{a_0}{2} + \sum_{k=1}^{\infty} (a_k \cos(kx) + b_k \sin(kx)),
 $$
then

$$
H_{\pi}[f](x) = \sum_{k=1}^{\infty} (a_k \sin(kx) - b_k \cos(kx)).
$$

Let us denote for arbitrary function $ \psi \in G\Psi(1,\infty) $ the new function
from at the same set
$$
\psi_H(p) = h(p) \cdot \psi(p).
$$

{\bf Theorem 11 }

$$
|| R_j \ f||G(\psi_H) \le Z \cdot ||f||G(\psi),
$$
where in general case the exact value of the constant $ Z $ is equal to $ 2/\pi: $

$$
\overline{Z} \stackrel{def}{=}
\sup_{\psi \in G\Psi(1,\infty)} \sup_{f \in G\psi(1,\infty), f \ne 0}
\frac{ ||R_j f||G(\psi_H)}{||f||G(\psi)} = \frac{2}{\pi}.
$$

{\bf Proof.} The upper bound for the constant $ Z $ it follows from the exact value
for the function $ h(\cdot); $ it remain to provide the low bound. \par
 It is sufficient to consider only the one-dimensional case, i.e.  when the Riesz
 transform coincides with the Hilbert transform. \par

{\bf A.} We consider at first the case $ X = (-\pi,\pi). $ Let us consider the
{\it family } of a functions of a view:

 $$
 g_{\Delta}(x) = \sum_{n=2}^{\infty} n^{-1} (\log n)^{\Delta} \sin(nx).
 $$
 Here $ \Delta = \const $ is some "great" constant. \par
   It is proved in \cite{Zygmund1}, p. 182-185 that  as $ x \to 0 $
 $$
 |g_{\Delta}(x)| \sim 0.5 \pi  |\log |x| \ |^{\Delta},
 $$
therefore

$$
| (2/\pi) \ g_{\Delta} |^p_p \asymp 2 \int_0^1 |\log x|^p dx = 2 \Gamma(\Delta p + 1);
$$
and as $ p \to \infty $
$$
| (2/\pi) \ g_{\Delta} |_p \sim \Gamma^{1/p}(\Delta p+1).
$$
 The correspondent Hilbert transform for the function $ g_{\Delta}, $ which we denote
 by $ - u_{\Delta}, $ has a view
 $$
 u_{\Delta}(x) = \sum_{n=2}^{\infty} n^{-1} (\log n)^{\Delta} \cos(nx).
 $$
 It is proved also in the book of Zygmund \cite{Zygmund1}, p. 182-185 that  as
 $ x \to 0 $

 $$
 |u_{\Delta}(x)| \sim (\Delta+1)^{-1} |\log |x| \ |^{\Delta+1},
 $$
 and following as $ p \to \infty $

 $$
 |u_{\Delta}|_p \sim (\Delta+1)^{-1} \Gamma^{1/p}((\Delta+1)p + 1).
 $$

 Let us denote
 $$
 V^{(\infty)}(g,p) = \frac{|T_H[g]|_p}{|h(p) \ g|_p}, \ p \in (1,\infty),
 $$
 $$
 \overline{V^{(\infty)}} = \sup_{g \in L(1,\infty) }
 \sup_{p \in (1,\infty)} V^{(\infty)}(g,p).
 $$
  It follows from the upper bounds of this theorem that
  $ \overline{V^{(\infty)}} \le 2/\pi.$
  On the other hand, we have for all the values $ \Delta = \const > 0 $ using the
  famous Stirling's formula:

 $$
 \overline{V^{(\infty)}} \ge \overline{\lim}_{\Delta \to \infty} \overline{\lim}_{p \to \infty} \frac{|T_H[g_{\Delta}]|_p}{|h(p) \ g_{\Delta}|_p} \ge
 $$

 $$
 \frac{2}{\pi} \overline{\lim}_{\Delta \to \infty} \overline{\lim}_{p \to \infty}
 \frac{\Gamma^{1/p}((\Delta+1)p + 1)}{p \cdot \Gamma^{1/p}(\Delta p + 1)} \ge
 $$

 $$
 \overline{\lim}_{\Delta \to \infty} \frac{2}{\pi}  e^{-1}
 \left( \frac{\Delta+1}{\Delta} \right)^{\Delta} =
\lim_{\Delta \to \infty} \frac{2}{\pi}  e^{-1}
 \left( \frac{\Delta+1}{\Delta} \right)^{\Delta} = \frac{2}{\pi}.
 $$
 Choosing as a function $ \psi(p) $ the following expression:

 $$
 \psi(p) := \psi_0(p) \stackrel{def}{=} |g_{\Delta}|_p, \ p \in (1,\infty),
 $$
for sufficiently  great values $ \Delta, $ we complete the consideration of the
case {\bf A} of this theorem. \par

\vspace{3mm}

{\bf B.} In this pilcrow we consider the case of the Hilbert's transform on the
whole real axis. Let us choose the function

$$
f_1(x) = I(x \in (c,d)), \ x \in R,  \ c,d = \const, \ d = c + 1;
$$
then $ \forall p \ge 1 \ |f_1|_p = 1. $ \par
 The correspondent Hilbert transform for the function $ f_1(x) $  we denote
 by $ v_1(x); $ it is equal to
 $$
 v_1(x) =  \frac{1}{\pi} \log \left| \frac{x-c}{x-d } \right|,
 $$
 see \cite{Bennet1}, p. 143-144. \par
 It is easy to see that as $ |x| \to \infty $
$$
v_1(x) \sim \frac{1}{\pi} \ |x|^{-1},
$$
and hence

$$
|v_1|_p \sim \frac{2}{\pi} \ \frac{1}{p-1}, \ p \to 1 + 0.
$$

 Let us denote as before
 $$
 V_1(f,p) = \frac{|T_H[f]|_p}{|h(p) \ f|_p}, \ p \in (1,\infty),
 $$
 $$
 \overline{V_1} = \sup_{f \in L(1,\infty) }
 \sup_{p \in (1,\infty)} V_1(f,p).
 $$
  It follows from the upper bounds of this theorem that  $ \overline{V_1} \le 2/\pi,$
  but
  $$
  \overline{V_1} \ge \overline{\lim}_{p \to 1+0} V_1(f_1, p) = \frac{2}{\pi}.
  $$
 Thus, $ \overline{V_1} = 2/\pi. $\par
 This completes the proof of our theorem.\\

 \hfill $\Box$ \\

\bigskip

\section{Discrete case}

\vspace{3mm}

We consider in this section the case when $ X = \{1,2,3, \ldots $ and $ \mu $ is  ordinary counting measure. So, the $ L_p $ norm of the function $ f = f(k), \
k \in X $ may be defined as follows:
$$
|f|_p = \left( \sum_{k=1}^{\infty} |f(k)|^p \right)^{1/p}, \ p \ge 1.
$$
In the terminology of a book  \cite{Bennet1} these spaces are called
{\it discrete resonant measurable spaces.} \par

 First of all we will formulate some simple properties of the $ L_p $ norms in
this space $ X. $ \par

$$
{\bf 1.} p \le q \ \Rightarrow |f|_q \le |f|_p;
$$

{\bf 2.} If for some $ p < \infty \ |f|_p < \infty, $ then

$$
\lim_{p \to \infty} |f|_p = \sup_{k \ge 1} |f(k)| \stackrel{def}{=} |f|_{\infty};
$$

{\bf 3.}  Natural functions. \\

Let $ f \in l_s $ for some value $ s, \ s \in [1, \infty); $ then $ \forall
q \ge p \ |f|_q \le |f|_s. $ Therefore, if we define the natural function for
the vector $ f, $ i.e. the function

$$
\psi_f(p) = |f|_p, \ s < p \le \infty,
$$
then the function  $  \psi_f(p)$ is bounded in the set $ p \in (s_1, \infty) $
and moreover

$$
\exists \lim_{p \to \infty} \psi_f(p) = \sup_{k} |f(k)| < \infty.
$$
{\it  We will consider in this section only the functions with the last both
Properties. } \par
{\bf Example.}  Let $ f(k) = k^{-\alpha}, \ \alpha \in (0,1). $ We get:

$$
\psi_f(p) \asymp \left[ \frac{ p }{p \ - \ 1/\alpha} \right]^{1/\alpha}, \ p \in
(1/\alpha,\infty).
$$

{\bf 4.} Tail behavior; see for comparison \cite{Ostrovsky2}. \par
 Let us denote for the (infinite) sequence $ f = \{f(k) \} $ the so-called tail
 function

$$
Z_{f}(\epsilon) = \mu \{k: |f(k)| > \epsilon \} =
 \card \{k: |f(k)| > \epsilon \}, \ \epsilon \ge 0.
$$
 Obviously, $ Z_f(\epsilon)= 0, $ if $ \epsilon > \sup_k |f(k)|. $ \par

 It follows from the Tchebyshev's inequality that if the sequence
$ f = \{f(k)\} $ belong to the space $ G(\psi), $ then

$$
Z_f(\epsilon) \le \inf_p
\left[ \frac{|f|_p^p \cdot \psi^p(p)}{\epsilon^p} \right].
$$

 Inversely,

$$
|f|_p =  \left[p \int_0^{\infty}y^{p-1} \ T(y) \ dy \right]^{1/p},
$$
Therefore

$$
||f||G(\psi) = \sup_{p \in (a,b)}
\left\{ \left[p \int_0^{\infty}y^{p-1} \ T(y) \ dy \right]^{1/p}
 /\psi(p) \right \}.
$$

{\bf Example.} If $ f(k) = k^{-\alpha}, k=1,2,\ldots, \alpha \in (0,1); $ then
$$
 \psi_f(p) \asymp \left[ p/(p-1/\alpha) \right]^{1/\alpha}; \
$$
$$
 Z_f(\epsilon) \asymp \epsilon^{-1/\alpha}, \ \epsilon \to 0+.
$$
But if it is given that
$$
|f|_p \le \left[ \frac{p}{p-1/\alpha} \right]^{1/\alpha},
$$
then we can conclude only

$$
Z_f(\epsilon) \asymp C \cdot \left[\frac{\epsilon}
{|\log \epsilon|}  \right]^{-1/\alpha}, \ \epsilon \in (0,1/e).
$$

\vspace{3mm}

 In the classical book \cite{Hardy1} there there are many examples of $ l_p \to l_p $
 norm estimations for matrix linear operators, which are completely analogous to the
 "continuous" case, as in section 3. For example:

$$
T^{(d)}_0[f](n) = n^{-1} \sum_{k=1}^{n} f(k), \ |T^{(d)}_0|_{p,p} = p/(p-1);
$$

 $$
 T^{(d)}_+[f](n):= \sum_{k=1}^{\infty} \frac{f(k)}{k+n}, \ |T^{(d)}_+|_{p,p} =
 \frac{\pi}{\sin(\pi/p)}, \ p \in (1,\infty) ;
 $$
 $$
 T^{(d)}_m[f](n):= \sum_{k=1}^{\infty} \frac{f(k)}{\max(k,n)}, \ |T^{(d)}_m|_{p,p} =
 \frac{p^2}{p-1}, \ p \in (1,\infty),
 $$
 $$
 T^{(d)}_d[f](n) := \sum_{k=n}^{\infty} \frac{f(k)}{k}, \ |T^{(d)}_d|_{p,p}=p, \ p \in (1,\infty),
 $$
$$
T^{(d)}_l[f](n) = \sum_{k=1}^{\infty} \frac{\log(k/n) f(k)}{k-n}, \ 0/0 = 0, \
|T^{(d)}_l|(p,p) = \left[ \frac{ \pi }{\sin(\pi/p) }  \right]^{2},
$$
  etc. \par

We will consider more generally linear matrix operators of a view:

$$
T^{(d)}_H[f](n) = n^{-1} \sum_{k=1}^{\infty} f(k) H(k/n),
$$
where the function $ H = H(z), z \in (0,\infty) $ satisfies the conditions of lemmas 1a
or 1b or 2a or 2b.  Suppose in addition the function $ H(\cdot) $ is piecewise
strictly monotonically decreasing  and piecewise continuous with finite points of
charging.\par

{\bf Theorem 12.}
$$
||T^{(d)}_H \ f||G(\psi_{(\phi)})  \le V_H \ ||f||G(\psi),
$$
where the exact value of constant $ V_H $ is equal to one.\par
{\bf Proof } is at the same as  in the proof of theorems 3 and 4, with at the same
"counter-examples"; more exactly, if the function $ f = f(x), x \in (0,\infty) $
is some example in the theorems 3 and 4, then the sequence
$$
f_0(k) = f(k)
$$
is the correspondent "counter-example" in the discrete case. \par

\vspace{3mm}

 Many another generalizations see in the works \cite{Bennet2}, \cite{Bennet3},
 \cite{Bennet4}, \cite{Bennet5}, \cite{Gao1}, \cite{Gao2}; see also reference
 therein. \par
  We will consider the linear operator with lower triangular (infinite) matrix
$ A = \{ a(n,k) \}$  with non-negative entries: $ a(n,k) \ge 0 $  of a view:

$$
a(n,k) = \lambda(k)/\Lambda(n), 1 \le k \le n; \ a(n,k) = 0, k \ge n+1;
$$
here
$$
\lambda(k) \ge 0, \lambda(1) > 0, \ \Lambda(n) = \sum_{k=1}^n \lambda(k).
$$
The correspondent linear operator $ T = T_A $ may be defined as ordinary, indeed,
for the infinite sequence $ x = \{ x(i), \ i = 1,2,3,\ldots  \} $

$$
T_A[x](j) = \sum_{i=1}^j a(j,i) x(i).
$$
It is proved in the article \cite{Bennet2}, see also \cite{Gao1}, \cite{Gao2}
that if

$$
L:=
\sup_n \left( \frac{\Lambda(n+1)}{\lambda(n+1)} - \frac{\Lambda(n)}{\lambda(n)}\right)
\in (1, \infty),
$$
then for the values $ p, \ p > L $

$$
|T_A|_{p,p} \le \frac{ p }{p-L}.
$$

As a consequence, if the function $ \psi(p) $ belongs to the set $ G\Psi(L,\infty) $
and if we define a new function

$$
\psi^{(L)}(p) = \frac{ p }{p-L} \ \psi(p),
$$
then it follows from theorem 1 that for arbitrary sequence
$  x = \{ x(i), i=1,2,\ldots   \} $ belonging to the space $ G(\psi) $

$$
||T_A[x]||G(\psi^{(L)}) \le 1 \cdot ||x||G(\psi). \eqno(33)
$$
{\bf Theorem 13}. The constant 1 in the inequality (33) is in general case exact.\par

\vspace{3mm}

 Analogously may be considered the case when $ \alpha = \const \in (0,1) $ and

 $$
 a(n,k) = n^{-\alpha} \left[k^{\alpha} - (k-1)^{\alpha} \right]
 $$
or

$$
a(n,k) = \frac{k^{\alpha-1}}{ \sum_{i=1}^n i^{\alpha-1}}.
$$
 In both the last cases the $ l_p \to l_p $ norm of the operator with the matrix
$  A = \{ a(n,k) \} $ allows the following estimation:

$$
|T_A|_{p,p} \le \frac{\alpha p}{\alpha p - 1}, \ p > 1/\alpha.
$$

As before, if the function $ \psi(p) $ belongs to the set $ G\Psi(1/\alpha,\infty) $
and if we define a new function

$$
\psi_{(\alpha)}(p) = \frac{ \alpha p }{\alpha p - 1} \ \psi(p),
$$
then it follows from theorem 1 that for arbitrary sequence
$  x = \{ x(i), i=1,2,\ldots   \} $ belonging to the space $ G(\psi) $

$$
||T_A[x]||G(\psi_{(\alpha)}) \le 1 \cdot ||x||G(\psi). \eqno(34)
$$
{\bf Theorem 14}. The constant 1 in the inequality (34) is exact.\par
{\bf Proof} of theorems 13 and 14. It is sufficient to prove the last theorem.
Let us consider as ordinary the following example:

$$
x = \{ x(k) \}, \ x(k) = k^{-\alpha} (\log k)^{\Delta}, \  k=1,2, \ldots; \
\Delta = \const > 0.
$$
 We get for the values $ p \to 1/\alpha + 0:$
$$
|x|^p_p \sim (\alpha p - 1)^{-1- \Delta p} \ \Gamma(\Delta p + 1),
$$

$$
 |x|_p \sim (\alpha p - 1)^{-\Delta -1/p} \ \Gamma^{1/p}(\Delta p + 1)
 \sim (\alpha p - 1)^{- \Delta - \alpha} \Gamma^{\alpha}(\Delta/\alpha + 1);
$$

$$
y := T_A[x] = \{ y(n) \}, \ n \to \infty \ \Rightarrow y(n) \sim  (\Delta + 1)^{-1}
\alpha \ n^{-\alpha} \ (\log(n))^{\Delta+1};
$$

$$
|y|_p \sim (\Delta+1)^{-1} \ \alpha \ (\alpha p - 1)^{-\Delta - \alpha - 1} \
\Gamma^{\alpha}( (\Delta+1)/\alpha + 1);
$$

$$
\overline{\lim}_{p \to 1/\alpha + 0}
 \left[  \frac{|y|_p}{|x|_p}: \frac{ \alpha p }{\alpha p - 1} \right]
\ge e^{-1} \left( 1 + \frac{1}{\Delta}  \right)^{\Delta};
$$
$$
\lim_{\Delta \to \infty} e^{-1} \left( 1 + \frac{1}{\Delta}  \right)^{\Delta} = 1.
$$

\hfill $\Box$ \\

\bigskip

 \section{Concluding remarks}

\vspace{3mm}

 The assertion of theorem 1 is true still for the so-called sublinear operators, i.e.
for the operators $ T[f] $  with properties:

$$
|T[f+g]|_p \le |T[f]|_p + |T[g]|_p, \ p \in (1,\infty);
$$

$$
|T[\lambda f ]|_p  = |\lambda| \ |T[f]|_p, \ \lambda = \const.
$$

 Many examples of sublinear operators give us maximal operators. Let
 $ T_{\theta}[\cdot], \ \theta \in \Theta $ be a family of linear operators. We put:

$$
T[f](x) = \sup_{\theta \in \Theta} | T_{\theta}[f](x) |,
$$
if there exists and satisfies, e.g. the inequality

$$
|T[f]|_p \le \Phi(p) \ |f|_p,
$$
where $ 0 < \Phi(p) < \infty, \ p \in (1,\infty). $\par
 For instance, let $ X = (-\pi, \pi) $  and let for any function $ f: X \to R \
 S^*[f] $ be the maximum of the absolute value of the partial Fourier $ S_n[f] $
 sum for the function $ f: $
 $$
 S^*[f] = \sup_n |S_n[f](x)|.
 $$
 It is known, see for example, \cite{Reyna1}, p. 151, that for some absolute constant
 $ C_1 \in (0,\infty) $

 $$
 |S^*[f]|_p \le C_1 \frac{p^4}{(p-1)^3} |f|_p, \ p \in (1,\infty),
 $$
 i.e. in this case $ \Phi(p) = p^4/(p-1)^3. $ \par
  Another example. Let $ H^*[f] $ be the maximal Hilbert transform for the function
  $ f, \ f: R \to R. $ It is proved in \cite{Reyna1}, p. 84 that for some absolute
  constant $ C_2 \in (0,\infty) $

 $$
 |H^*[f]|_p \le C_2 \frac{p^2}{p-1}|f|_p,
 $$
 i.e. here $ \Phi(p)= p^2/(p-1). $ \par

 It follows from theorem 1 that if in the case when $ \Phi(p) < \infty, \ p \in
 (1,\infty) $ and $ \psi(\cdot) \in G\Psi(1,\infty) $ we define

 $$
 \tilde{\psi}(p) = \Phi(p) \cdot \psi(p),
 $$
 then

 $$
 ||T[f]||G( \tilde{\psi}) \le ||f||G(\psi).
 $$

 Analogous asymptotically as $ p \to 1+0 $ and $ p \to \infty $ exact estimations
for the function $ \Phi = \Phi(p) $ are obtained  for  the famous singular operators
of Calderon - Zygmund type or analogously maximal Calderon - Zygmund type,
see  \cite{Adams1},  \cite{Bennet1}, \cite{Dunford1}, \cite{Stein1} etc. \par

\end{document}